\documentclass[12pt]{amsart}
\usepackage{amssymb,array}

\newtheorem{lemma}{Lemma}
\newtheorem{proposition}[lemma]{Proposition}
\newtheorem{theorem}[lemma]{Theorem}

\newcommand{\CE}{\mathcal{E}}
\newcommand{\CH}{\mathcal{H}}
\newcommand{\CL}{\mathcal{L}}

\newcommand{\R}{\mathbb{R}}
\newcommand{\C}{\mathbb{C}}

\newcommand{\N}{\mathbb{N}}

\newcommand{\T}{\mathbb{T}}

\newcommand{\spec}{\mathrm{Spec}\,}
\newcommand{\spectwo}{\mathrm{Spec}_2\,}
\newcommand{\num}{\mathrm{Num}\,}
\newcommand{\dom}{\mathrm{Dom}\,}
\newcommand{\ran}{\mathrm{Ran}\,}
\newcommand{\spa}{\mathrm{Span}\,}
\newcommand{\conv}{\mathrm{Conv}\,}
\renewcommand{\Re}{\mathrm{Re}\,}
\renewcommand{\Im}{\mathrm{Im}\,}
\newcommand{\dist}{\mathrm{dist}\,}

\title[]{Limiting set of second
order spectra}
\author[L. Boulton]{Lyonell Boulton}

\date{$1^{th}$ August 2003}



\begin{document}
\begin{abstract}
Let $M$ be a self-adjoint operator acting on a Hilbert space
$\CH$. A complex number $z$  is in the second order spectrum of
$M$ relative to a finite dimensional subspace $\CL\subset \dom
M^2$, iff the truncation to $\CL$ of $(M-z)^2$ is not invertible.
This definition was first introduced in \cite{secr} and according
to the results of \cite{lesh} and \cite{shar}, these set provide a
method for estimating eigenvalues free from the problems of
spectral pollution. 
In this paper 
we investigate various aspects related to the issue of approximation 
using second order spectra.  Our main result shows 
that under fairly mild hypothesis on $M,$ the uniform
limit of these set, as $\CL$ increases towards $\CH$, contains the
isolated eigenvalues of $M$ of finite multiplicity. Therefore, unlike the
majority of the standard methods, second order spectra combine
non-pollution and approximation at a very high level of
generality.
\end{abstract}

\maketitle

\section{Introduction}

Let $M$ be a self-adjoint operator acting on an infinite
dimensional Hilbert space $\CH$ and let $\lambda$ be an isolated
eigenvalue of $M$ of finite multiplicity. Suppose we intend to
estimate $\lambda$ numerically. A natural method, which is already
over one hundred years old, is to truncate $M$ by taking an
orthonormal basis $\{\phi_j\}_{j=1}^\infty$ of $\CH$ contained in
the (operator or form) domain of $M$ and compute the eigenvalues
of the $n\times n$ matrix $M_n=(\langle M\phi_j,\phi_k
\rangle)_{j,k=1}^n$ for large $n$, expecting that some of these
eigenvalues would be close to $\lambda$. Unfortunately, an ill
chosen $\{\phi_j\}_{j=1}^\infty$ will result in either of the
following difficulties:
\begin{itemize}
\item[-] \emph{Lack of approximation:} no
eigenvalue of $M_n$ is close to $\lambda$. This happens when $M$
is unbounded. \pagebreak
\item[-] \emph{Spurious eigenvalues:} or
spectral pollution, $M_n$ has many eigenvalues unrelated to the
spectrum of $M$ in a large neighbourhood of $\lambda$. This occurs
when $\lambda$ is in a gap, i.e. in between two points, of the
essential spectrum of $M$.
\end{itemize}

It is well known that if $M$ is an operator bounded below and
\linebreak $\lambda < \min [\spec\!_{\mathrm{ess}}\, M]$, then the
Rayleigh-Ritz theorem provides satisfactory procedures to deal
with these two issues. For instance, suppose that $\lambda$ is the
first eigenvalue of $M$ and it is non-degenerate. If
$\{\phi_j\}_{j=1}^\infty$ is such that
\begin{equation} \label{e1}
M_n\left(\sum_{j=1}^n \langle\psi, \phi_j\rangle \phi_j \right)
\to \lambda \psi \ \mbox{as}\ n\to \infty\ \mbox{whenever}\
M\psi=\lambda \psi,
\end{equation}
then the first eigenvalue of $M_n$ converges to $\lambda$ (i.e. we
achieve approximation) and the second eigenvalue of $M_n$, counting
multiplicity, can not be smaller than the second eigenvalue of $M$
(i.e. no chance of spurious eigenvalues). Condition \eqref{e1} is
useful in applications, because it can sometimes be verified on
abstract grounds without having much information about the
eigenfunctions $\psi$, cf. \cite[theorem XIII.4]{resi}.

In contrast, when $\lambda$ is in between two points of the
essential spectrum, it is well known that pollution can arise. In
fact, for each $\mu\in \R\setminus \spec M$ lying between two
points of the essential spectrum, there exists an orthonormal
basis $\{\phi_j\}_{j=1}^\infty$ and a subsequence
$n(1)<n(2)<\ldots$, such that $\mu \in \spec M_{n(k)}$ for all
$k\in\N$ (see \cite{lesh}). There is a large amount of literature
devoted to study this problem in applications, we refer to
\cite{sp} and \cite{lesh}, and the references therein, for a more
complete account on the matter.

In order to find eigenvalues in gaps of the essential spectrum,
one can appeal to the following geometrically motivated idea (see
\cite{sp} and \cite{kato2}). Suppose $a<\zeta<b$, such that
$a,b\!\in\! \spec_{\mathrm{ess}}M$ but $(a,b)\cap\spec_{\mathrm{ess}}
M=\varnothing$. Then $(M-\zeta)^2$ is bounded below by zero and
the eigenvalues below $\spec_{\mathrm{ess}}(M-\zeta)^2$ are the
$(\lambda -\zeta)^2\geq 0$, where $\lambda$ is an eigenvalue of
$M$ inside $(a,b)$. Hence by applying Rayleigh-Ritz principle to
\linebreak $(M-\zeta)^2$, we might be able to estimate (up to a
square root ambiguity) the eigenvalue $\lambda$. If, for instance,
$(M-\zeta)^2$ satisfies \eqref{e1}, from the first eigenvalue of
its matrix truncations we may estimate the closest eigenvalue of
$M$ to $\zeta$ in $(a,b)$.

Instead of looking for the eigenvalues close to zero for a
truncation of $(M-\zeta)^2$, we can also consider finding $z$ such
that the truncation of $(M-z)^2$ is singular. Denote by $\CL_n$
the $\spa \{\phi_j\}_{j=1}^n$ and by $P_n$ the orthogonal
projection onto $\CL_n$. Assume that $\CL_n\subset \dom M^2$. We
say that $z\in \C$ is an element of the second order spectrum of
$M$ \linebreak relative to $\CL_n$, denoted below by $\spectwo
(M|\CL_n)$, if and only if \linebreak $\det [P_n(M-z)^2|\CL_n]=0$.
This definition was first introduced in \cite{secr} by Davies.
Levitin and Shargorodsky have recently proposed, see \cite{lesh},
using second order spectra to approximate isolated eigenvalues
inside spectral gaps, on the grounds that they are free from
pollution. To be precise, if $z\in \spectwo (M|\CL_n)$, then
\begin{equation} \label{e2}
   [\Re z-|\Im z|,\Re z+|\Im z|]\cap \spec M \not =\varnothing
\end{equation}
(cf. \cite{lesh} and \cite{shar}). In other words, a point in
$\spectwo (M|\CL_n)$ is close to the real axis only when it is
also close to the spectrum. Therefore numerical procedures based
on second order spectrum, never lead to spurious spectrum.

Although the bound \eqref{e2} does not guarantee that any point at
all in the spectrum is estimated for large $n$, some numerical
experiments performed in \cite{lesh} and the results of
\cite{bou}, indicate that approximation occurs in test models.

In this paper we discuss rigorously various aspects of the
question of approximation in the method considered by Levitin and
Shargorodsky. For this, we study the uniform limiting set
\[
\Lambda:=\mathrm{u}\!\!-\!\!\lim_{n\to\infty} \spectwo (M|\CL_n)=
\{z\in \C\,:\,z_n\to z,\, z_n\in \spectwo (M|\CL_n)\}.
\]
The interest in $\Lambda$ is clear, the points in $\spec M$ which
are estimated by $\spectwo(M|\CL_n)$ for $n$ large, are those in
$\Lambda\cap \R$. The following assumption on $\lambda$ and $P_n$
will be crucial in our subsequent analysis\vspace{.1in},
\newline \vspace{.1in} (H) \hspace{.5in} \begin{minipage}{3.2in}
\emph{if $M\psi=\lambda \psi$, \ then} \quad $\|P_nMP_n\psi -
\lambda \psi\|\to 0$ \quad \emph{and} \quad $\|P_nM^2P_n\psi -
\lambda^2 \psi\|\to 0,$ \ \emph{as \ $n\to \infty$.}\end{minipage}
\newline Our principal task will be to show the following

\begin{theorem} \label{t4}
If $\lambda$ is an isolated eigenvalue of finite multiplicity of
$M$ and (H) holds, then $\lambda \in \Lambda$.
\end{theorem}

Obviously any bounded operator satisfies (H), and notice the
analogy between (H) and \eqref{e1} in the general case. We will
discuss in section~5 how to verify this condition for unbounded
operators without having much information about $\psi$.

We should emphasize the relevance of this result. In many
applications, the essential spectrum can be found analytically,
whereas the isolated eigenvalues of finite multiplicity are the
ones that should be estimated numerically. Theorem~\ref{t4} shows
that second order relative spectra provide a method which will
give convergence to these eigenvalues. Furthermore, in conjunction
with \eqref{e2}, it is ensured that the method will not pollute,
no matter the location of these eigenvalues. We are not aware of
any other technique that combines approximation and non-pollution
at such level of generality.

In section~2 we keep close to the ideas of \cite{secr}. We
characterize $\spectwo(M|\CL_n)$ and the set $\Lambda$ in terms of
resolvent norms. This characterization will be important in
section~3, where we prove theorem \ref{t4}. In section~4 we
discuss in detail a simple model bounded operator, for which
$\Lambda$ can be found explicitly. In order to test
theorem~\ref{t4} in concrete situations, we also include some
experimental numerical outputs for rank one perturbations of this
model. Finally section~5 is devoted to discuss the verification of
condition (H) for concrete unbounded situations.


\section{Second order spectrum and norm of inverses}

The sequence of subspaces $\CL_n$ does not necessarily need to be
nested in order to accomplish approximation. Below $P_n$ denotes a
family of orthogonal projections and $\CL_n:=\ran P_n\subset \dom
M^2$, such that $\dim \CL_n=n$ and $P_nx\to x$ for all $x\in \CH$.
For the linear operator $H$, we shall often denote by $H_n$ the
truncation $P_nH|\CL_n$. By $H_n^{-1}$ we shall always mean the
inverse of $H_n$ within $\CL_n$. If $H_n$ is not invertible we
allow the abuse of notation $\|H_n^{-1}\|=\infty$.

From the definition it is easy to deduce that $\spectwo(M|\CL_n)$
is a set of at most $2n$ different complex numbers. In general
these points do not intersect the real line, unless $\CL_n$
contains an eigenfunction of $M$. Since
\[\overline{\det P_n(M-z)^2|\CL_n} = \det P_n(M-\overline{z})^2|\CL_n,\]
$\spectwo(M|\CL_n)$ is symmetric with respect to $\R$. The
limiting set $\Lambda$ might or might not contain non-real points.
The reference \cite{bou} is largely devoted to showing that
$\Lambda\subseteq \spec M$, when $M$ is a discrete Schr\"odinger
operator whose potential decays fast at infinity. On the other
hand, in the example we discuss in section~4, $\Lambda$ is equal
to the unit circle.

In order to describe properties of $\Lambda$, it is useful
characterize the points in $\spectwo(M|\CL_n)$ as the zeros of a
certain Lipschitz continuous function $\sigma_n(z)$. Davies in
\cite{secr}, and Davies and Plum in \cite{sp} study a procedure
for computing eigenvalues of bounded self-adjoint operators, based
on evaluating $\sigma_n(z)$ at $z\in \R$. Let
\[
   \sigma_n(z):=\inf \left\{ \frac{\|P_n(M-z)^2|\CL_n
   v\|}{\|v\|}\,:\, 0\not = v \in \CL_n \right\}.
\]
If $z\not \in \spectwo(M|\CL_n)$, then
$\sigma_n(z)^{-1}=\|[P_n(M-z)^2|\CL_n]^{-1}\|$. Moreover,
$\sigma_n(z)=0$ if and only if $z\in \spectwo (M|\CL_n)$. Both
assertions are direct consequence of the definition.

The following result shows that, in order to find the zeros of
$\sigma_n$, we can start from an initial guess and move in steps
towards a local minimum.

\begin{lemma} \label{t3}
The function $\sigma_n(z)$ is a non-negative Lipschitz continuous
function of the variable $z$. Its only local minima are the points
where it vanishes.
\end{lemma}
\proof We follow closely \cite[theorem 21]{secr}. If
$|w-z|<\varepsilon$, then
\begin{equation} \label{e12}
\begin{aligned}
  \sigma&_{n}(z) =\inf \frac{\|P_{n}(M-z)^2|\CL_{n}v\|}
  {\|v\|} \\
  & \leq \inf \frac{\|P_{n}(M-w)^2|\CL_{n}v\|+|w-z|\,
  \|P_{n}(2M-z-w)|\CL_{n}v\|}{\|v\|} \\
  & \leq \sigma_{n}(w)+|w-z|\sup \frac{\|(2M_{n}-z-w)v\|}
  {\|v\|} \\
   &\leq \sigma_{n}(w)+|w-z|c,
\end{aligned}
\end{equation}
where $c>0$ can be chosen independent of $\varepsilon$. This
ensures the Lipschitz continuity of $\sigma_n$.

On the other hand, given $z_0 \in \C \setminus \spectwo(M|\CL_n)$,
choose $u,v\in \CL_n$ vectors of norm 1 such that
\[
   \sigma_n(z_0)^{-1}=|\langle [P_n(M-z_0)^2|\CL_n]^{-1}u,v\rangle|.
\]
The function $\langle [P_n(M-z)^2|\CL_n]^{-1}u,v\rangle$ is
analytic in $z \in \C \setminus \spectwo(M|\CL_n)$. We claim that
it is not constant. If such a claim is true, then there always
exists $z$ in any neighbourhood of $z_0$ such that
\begin{align*}
    \sigma_n(z)^{-1} & \geq \big|\langle [P_n(M-z)^2|
    \CL_n]^{-1}u,v\rangle \big| \\
    & > \big|\langle [P_n(M-z_0)^2|
    \CL_n]^{-1}u,v\rangle \big| \\
    & = \sigma_n(z_0)^{-1},
\end{align*}
ensuring the second assertion of the lemma. In order to prove
our claim, first notice that
\begin{align*}
    \big|\langle [P_n(M&-z)^2|\CL_n]^{-1}u,v\rangle \big|  \leq
   \|[P_n(M-z)^2|\CL_n]^{-1}\| \\
   & = |z|^{-2}\|[z^{-2}P_nM^2|\CL_n-2z^{-1}P_nM|\CL_n+1]^{-1}\|.
\end{align*}
Then, if
$
   |z|>3 \max\{\|P_nM^2|\CL_n\|,\|P_nM|\CL_n\|,1\},
$
\begin{align*}
  \|z^{-2}P_nM^2|\CL_n-2z^{-1}P_nM|\CL_n\|<7/9<1
\end{align*}
and hence,
\[
   \big|\langle [P_n(M-z)^2|\CL_n]^{-1}u,v\rangle \big|\leq
   (9/2)|z|^{-2} \to 0 
\]
as $|z|\to \infty$. Thus necessarily the abovementioned function is not 
constant. 
\medskip

In conjunction with $\Lambda$, it is natural to consider the set
of asymptotic zeros of the family $\{\sigma_n\}$,
\[
  \Sigma\equiv \Sigma(M):=\{z\in \C:\sigma_n(z)\to 0,\,n\to \infty\}.
\]
If $M$ is bounded, this set contains the limiting set $\Lambda$.
Indeed, suppose that $z\not \in \Sigma$. Then there is a
subsequence $n(j)$ and $c>0$, such that
$P_{n(j)}(M-z)^2|\CL_{n(j)}$ are invertible for all $j\in \N$ and
\[
c\geq \sigma_{n(j)}^{-1}(z)=\|[P_{n(j)}(M-z)^2|\CL_{n(j)}]^{-1}\|.
\]
By virtue of \eqref{e12}, for $|z-w|<\varepsilon$,
\[
    \sigma_{n(j)}(z)\leq \sigma_{n(j)}(w)+|w-z|c_1
\]
where $c_1>0$ can be chosen independent of $\varepsilon$ and $j$
(here we use that $M$ is bounded in order to ensure that $c_1$
does not depend on $j$). By letting $\varepsilon$ be small enough,
we can find $c_2>0$ such that $\sigma_{n(j)}(w)\geq c_2$ whenever
$|w-z|<\varepsilon$. Thus
\[
    \bigcup_{j=1}^\infty \spectwo(M|\CL_{n(j)})\cap
    \{|z-w|<\varepsilon\} =\varnothing
\]
so therefore $z\not \in \Lambda$.

Furthermore, whenever $M$ is bounded, $\Sigma\,\cap\, \R=\spec M$.
Indeed, if $\lambda\in \spec M$, for each $k>0$ there is
$\psi_k\in \CH$, $\|\psi_k\|=1$, such that
$\|(M-\lambda)^2\psi_k\|<1/k$. Then
\[
   \lim_{n\to\infty}
   \frac{\|P_n(M-\lambda)^2P_n\psi_k\|}{\|P_n\psi_k\|}< 1/k
\]
and so $\sigma_n(\lambda)\to 0$ as $n\to\infty$. Conversely,
notice that
\begin{equation} \label{e6}
   \left<P_n(M-z)^2|\CL_nv,v\right>=\left<(M-z)^2v,v\right>
\end{equation}
for all $v\in \CL_n$, $\|v\|=1$. Then, if $\lambda\in \R$ but
$\lambda\not \in \spec M$,
\[
  \mathrm{Num}\, P_n(M-\lambda)^2|\CL_n \subset
  \mathrm{Num}\, (M-\lambda)^2 \subset [b,\infty) 
\]
for some constant $b>0$ independent of $n$. Hence 
\[
  \sigma_n(\lambda) \geq \dist[0,\mathrm{Num}\,P_n(M-\lambda)^2|\CL_n]
   \geq b,
\]
so that $\lambda \not \in \Sigma \cap \R$. This observation is
crucial in the method of Davies and Plum mentioned earlier. Notice
that the direction $\Sigma\,\cap\, \R\subseteq \spec M$ does not
require $M$ to be bounded. Here and below ``Num'' denotes numerical range.

\medskip

It would be of obvious interest to find hypothesis which guarantee
$\Lambda=\Sigma$. We show that the equality holds at least in
particular cases. Below and elsewhere we denote by $\T$ the unit
circle.

\begin{proposition} \label{t1}
If $\spec M$ consists of two points $a<b$, then
\[\Sigma=\Lambda\subseteq \{z\in \C:|z-(a+b)/2|=(b-a)/2\}.\]
\end{proposition}
\proof Since
\[
   \spectwo(\alpha M+\beta|\CL_n)=\alpha \spectwo(M|\CL_n)+\beta
\]
and $\Sigma(\alpha M+\beta)=\alpha \Sigma(M)+\beta$ for scalars
$\alpha$ and $\beta$, it is enough to show the desired property
for the case $\spec M=\{\pm 1\}$. By virtue of the hypothesis, it
is clear that $z=0$ can neither be in $\Sigma$ nor in any of the
$\spectwo(M|\CL_n)$, so below we assume $z\not=0$.

The latter assumption on $M$ yields $M^2=\mathrm{Id}$. Then, since
\begin{equation} \label{e8}
   P_n(M-z)^2|\CL_n=2z\big[(z^{-1}+z)/2-P_nM|\CL_n\big],
\end{equation}
$z\in \spectwo (M|\CL_n)$ if and only if \[(z^{-1}+z)/2\in \spec
M_n\subset \mathrm{Num}\, M=[-1,1].\] Hence all the second order
relative spectra of $M$ are contained in $\T$. This shows the
second relation.

We already saw that $\Lambda \subseteq \Sigma$ in general. In
order to show the reverse inclusion, assume that $z\not \in
\Lambda$. Then there is a subsequence $n(j)$ and $\delta>0$, such
that
\[
   \bigcup_{j=1}^\infty \spectwo(M|\CL_{n(j)}) \cap \{|z-w|<\delta\}=
   \varnothing.
\]
By virtue of \eqref{e8},
\[
   \sigma_{n(j)}(z)=2|z|\,\dist \big[(z^{-1}+z)/2,\spec
   M_{n(j)}\big].
\]
But if $(w^{-1}+w)/2$ is in the spectrum of $M_{n(j)}$,
\begin{align*}
  |(z^{-1}-z)/2-(w^{-1}-w)/2|&=|z|^{-1}|w-z|\,|w^{-1}-z|\\ &\geq
  |z|^{-1}\delta^2.
\end{align*}
Hence $\sigma_{n(j)}(z)\geq c>0$ for all $j\in \N$ and so $z\not
\in \Sigma$. This shows proposition~\ref{t1}.

\medskip

In applications it is of interest to understand how the second
order spectrum and its limiting set change under compact
perturbations. The following lemma provides some lights in this
respect. Let $H$ be a bounded not necessarily self-adjoint
operator and $K$ be compact operator. It is well known that
$\|H_n^{-1}\|$ is uniformly bounded for all $n$ large and $H+K$ is
invertible, if and only if $\|(H+K)_n^{-1}\|$ is uniformly bounded
for all $n$ large and $H$ is invertible (see, for instance
\cite[theorem 2.16]{bosi}). This allows us to prove that $\Sigma$
does not change substantially under compact perturbation.

\begin{lemma} \label{t7}
If $M$ is bounded and $K=K^\ast$ is a compact operator, then
\[
   \Sigma(M)\cup \spec_{\mathrm{disc}}(M+K)=\Sigma(M+K)\cup
   \spec_{\mathrm{disc}}M.
\]
In particular $\Sigma(M)$ and $\Sigma(M+K)$ coincide outside the
real axis.
\end{lemma}
\proof The point $z$ is not in the set at the left hand side if
and only if,
\[
(M+K-z)^2=(M-z)^2+[(M-z)K+K(M-z)+K^2]
\]
is invertible and
\[
   \|[P_{n(j)}(M-z)^2|\CL_{n(j)}]^{-1}\|
\]
is uniformly bounded for some suitable subsequence $n(j)$.
According to the above observation, the latter is equivalent to
having $(M-z)^2$ invertible and
\[
   \|[P_{n(j)}(M+K-z)^2|\CL_{n(j)}]^{-1}\|
\]
uniformly bounded, so therefore it is equivalent to $z$ not being
in the set at the right hand side.


\section{Proof of theorem 1}

We start by noticing that condition (H) ensures that
\begin{equation} \label{e7}
\sigma_n(\lambda)\to 0\qquad \mathrm{as}\quad n\to \infty.
\end{equation}
Indeed, take $v=P_n\psi$, then
\begin{align*}
\sigma_n(\lambda) & \leq \frac{\|P_n(M-\lambda)^2P_n\psi\|}{\|P_n
\psi\|} \\&= \frac{\|P_nM^2P_n\psi-2\lambda P_nMP_n\psi+ \lambda^2
\psi\|}{\|P_n \psi\|} \to 0
\end{align*}
as $n\to \infty$.

Decompose
\begin{gather*}
   \CH=\CE \oplus \CE^{\perp} \\
   \CE= \spa\{\psi\in \dom M\,:\, M\psi=\lambda \psi \}
   \subset \dom M^2
\end{gather*}
where $\dim \CE<\infty$, and both $\CE$ and $\CE^\perp$ are
invariant under $M$. According to the standard notion, $\CE^\perp$
invariant under $M$ means that $Mx\in \CE^\perp$ for all $x\in
\CE^\perp \cap \dom M$. The restriction $M|\CE$ corresponds to
multiplication by $\lambda$. Denote by $P_\CE$ the orthogonal
projection onto $\CE$. Put
\[
    M=M(I-P_\CE)+(\lambda-\mu)P_\CE+\mu P_\CE=:A+\mu P_\CE
\]
where $0\not = \mu\not \in \spec M$ and $\dom A=\dom M$. Then
$A=A^\ast$, $M$ is a finite rank perturbation of $A$ and
\[
   \spec A=\{\lambda-\mu\} \cup \spec M \setminus \{\lambda\} .
\]

For $z\in \C$, let $A(z):=(A-z)^2$ with $\dom A(z)=\dom M^2$. Then
$A(z)$ is a holomorphic family of type~A for $z\in \C$. If
$K(z)=\mu(2\lambda-2z-\mu)P_\CE$, then
\[
   (M-z)^2x=A(z)x+K(z)x,\qquad \qquad x\in \dom M^2.
\]
Put $A_n(z):=P_nA(z)|\CL_n$ and $K_n(z)=P_nK(z)|\CL_n$.

\medskip

The validity of theorem~\ref{t4} can be formally justified by the
following observation. For $z$ close to $\lambda$, $\spec A(z)$ is
the spectrum of $A$ ``bent'' to the right half plane by the map
$w\mapsto (w-z)^2$. Since the numerical ranges of $A_n(z)$ are far
from zero (cf. lemma \ref{t5}), $\|A_n(z)^{-1}\|$ are uniformly
bounded for all such $z$. When we add $K(z)$, the truncation are
not sectorial any longer, but since $K(z)$ is of finite rank, for
$n$ large, $A_n(z)+K_n(z)$ are still invertible, except for a few
isolated points which correspond to perturbations of the isolated
eigenvalue $\lambda$ of $A(z)+K(z)$. When some of these isolated
points is equal to zero, $z\in \spectwo(M|\CL_n)$.

\begin{lemma} \label{t5}
Let $\delta>0$ be small enough. Then there exist non-negative
constants $\alpha,\, \beta$ and $b$, such that
\[
 \num A_n(z) \subseteq \conv \left[([\alpha,\infty)+
 i[-\beta,\beta])^2\right] \subset \{\Re(z)\geq b>0\}
\]
for all $n\in \N$ and $z\in \{|z-\lambda|\leq \delta\}$.
\end{lemma}
\proof By substituting $A$ for $M$ in \eqref{e6}, clearly  $\num
A_n(z) \subseteq \num A(z)$ for all $n\in \N$. Since $A=A^\ast$,
then $A(z)$ is a normal operator for all $z\in \C$, so that
\[
   \num A(z) \subseteq \overline{\conv \left[\spec (A-z)^2
   \right]}
   \qquad \qquad z\in \C.
\]
Since $\lambda \not \in \spec A$, we can find $\delta>0$ small
enough, such that $\{|z-\lambda|\leq \delta\}$ does not intersect
$\spec A$. Then there are $\alpha,\, \beta>0$, such that
\[
  \spec(A-z)\subset \big((-\infty,-\alpha]+i[-\beta,\beta]\big)\cup
  \big([\alpha,\infty)+i[-\beta,\beta]\big)
\]
for $|z-\lambda|\leq \delta$. Hence
\[
  \spec(A-z)^2 \subset ([\alpha,\infty)+i[-\beta,\beta])^2.
\]
The $b>0$ can be found by choosing $\delta>0$ small enough, such
that $\beta$ is also small.

\medskip

\begin{lemma} \label{t2}
Let $0\not = \psi \in \dom H$ be such that $H\psi=\nu \psi$ and
$H_nP_n \psi \to \nu \psi$. If $H$ and $H_n$ are invertible, and
$\|H_n^{-1}\|\leq c$ where $c>0$ is independent of $n$, then
\[
   \|H_n^{-1}P_n\psi - H^{-1}\psi\| \to 0
\]
as $n\to \infty$.
\end{lemma}
\proof
\[
   \|H_n^{-1}P_n \psi-H^{-1}\psi \| \leq c \|P_n \psi - \nu^{-1}H_n P_n
   \psi \| \to 0
\]
as $n\to \infty$.

\medskip

It is well known that for any bounded operator $H$ and $w\not \in
\num H$,
\[
   \|(H-w)^{-1}\|\leq \dist \left(w,\num H \right)^{-1}
\]
(cf. for instance \cite[p.268]{kato}). Then according to lemma
\ref{t5}, for all  $n\in \N$ and $|z-\lambda|<\delta$, $A_n(z)$
are invertible and there exists $c_1>0$ independent of $n$ and
$z$, such that $\|A_n(z)^{-1}\|<c_1$. Since $\CE$ is an eigenspace
of $A(z)$, (H) ensures that $A_n(z)P_n\psi\to A(z)\psi$ as $n\to
\infty$ for all $\psi\in \CE$. Since $z\not \in \spec A$, $A(z)$
is invertible. Therefore lemma \ref{t2} yields
\begin{equation} \label{e3}
   \|A_n(z)^{-1}P_n \psi - A(z)^{-1}\psi\|\to 0, \qquad \psi
   \in \CE
\end{equation}
for all $|z-\lambda|<\delta$.

\medskip

The various constants $c_j>0$ that appear below are independent of
$z$ and $n$.

In order to prove theorem \ref{t4}, we show that for $\delta>0$
small enough, we can find constants $c>0$ and $\tilde{N}>0$
uniform in $z$, such that $\sigma_n(z)^{-1}\leq c$ for all
$|z-\lambda|=\delta$ and $n\geq \tilde{N}$. In virtue of
\eqref{e7} and by applying lemma~\ref{t3}, the above ensures the
existence of $N>0$ such that we can find
\[z_n\in \{|z-\lambda|<\delta\}\cap \spectwo(M|\CL_n)\] for all
$n\geq N$.

Let $\delta>0$ be small enough such that lemma \ref{t5} holds and
\[
   \spec M \cap \{|z-\lambda|\leq \delta \} =\{\lambda\}.
\]
Denote by
\[
  N_\delta:=\{z\in \C\,:\,|z-\lambda|=\delta\}.
\]
Since $A(z)+K(z)=(M-z)^2$ is invertible, for all $z\in N_\delta$,
\[
   I+A(z)^{-1}K(z)=A(z)^{-1}(A(z)+K(z))
\]
is also invertible. Since  the function
$\|I+A(\cdot)^{-1}K(\cdot)\|:N_\delta\longrightarrow (0,\infty)$
is continuous, it should be uniformly bounded. Put $c_2>0$ such
that
\begin{equation} \label{e5}
  \|(I+A(z)K(z))^{-1}\| \leq c_2^{-1},\qquad \qquad z \in
  N_\delta.
\end{equation}

Let $\varepsilon>0$ and $z\in N_\delta$. Since $N_\delta$ is
compact, there exists a finite set of points $\{w_j\}\subset
N_\delta$ satisfying the following property. Given any $z\in
N_\delta$, there is $w\in \{w_j\}$ such that
\begin{equation} \label{e4}
   \|A(z)^{-1}-A(w)^{-1}\|<\varepsilon \qquad
   \mbox{and} \qquad |w-z|<\varepsilon.
\end{equation}
Since $\dim \CE<\infty$ and $\{w_j\}$ is finite, in virtue of
\eqref{e3}, there exists $\tilde{N}>0$ such that
\[
   \|A_n(w)^{-1}P_n\psi - A(w)^{-1}\psi\|<\varepsilon,
\]
for all $n\geq \tilde{N}$, $w\in (w_j)$ and $\psi\in \CE$. Thus,
if $w$ satisfies \eqref{e4},
\begin{align*}
   \|A_n&(z)^{-1}P_n\psi -A(z)^{-1} \psi\|
    \\ & \leq
   \|A_n(z)^{-1}P_n\psi -A_n(w)^{-1} P_n\psi\|+
    \|A_n(w)^{-1}P_n\psi -A(w)^{-1} \psi\|+ \\
   & \quad + \|A(w)^{-1}\psi -A(z)^{-1} \psi\| \\
   & < \|A_n(z)^{-1}P_n\psi -A_n(w)^{-1} P_n\psi\|+2\varepsilon \\
   &=\|A_n(z)^{-1}A_n(w)^{-1}(A_n(w)-A_n(z))P_n\psi\|+2\varepsilon  \\
  &\leq c_1^2\|(w-z)P_n[(w+z)-2A]P_n\psi\|+2\varepsilon \\
   & \leq |w-z|c_1^2c_3+2\varepsilon <c_4 \varepsilon.
\end{align*}
The existence of $c_3$ is ensured by (H). By choosing
$\varepsilon$ small enough, we can find $\tilde{N}>0$ independent
of $z$,  such that
\begin{equation*}
  \|A_n(z)^{-1}P_nK(z)-A(z)^{-1}K(z)\|<c_2/2, \qquad n\geq \tilde
  N.
\end{equation*}

Let $x\in \CH$. In virtue of \eqref{e5},
\begin{align*}
   c_2\|P_nx\|&\leq \|(I+A(z)^{-1}K(z))P_nx\| \\
   & \leq \|(I+A_n(z)^{-1}P_nK(z))P_nx\|+ \\
   & \quad +\|A_n(z)^{-1}P_nK(z)-A(z)^{-1}K(z)\|\|P_nx\|.
\end{align*}
Hence, for all $n\geq \tilde{N}$,
\begin{align*}
  (c_2/2)\|P_nx\| & \leq \|(I+A_n(z)^{-1}P_nK(z))P_nx\| \\
    & \leq \|A_n(z)^{-1}\| \|(A_n(z)+P_nK(z)P_n)P_nx\| \\
    &\leq c_1 \|P_n(M-z)^2|\CL_nP_nx\|.
\end{align*}
Therefore $P_n(M-z)^2|\CL_n$ is invertible and
\[
  \sigma_n(z)^{-1}=\|(P_n(M-z)^2|\CL_n)^{-1}\| \leq 2c_1/c_2
\]
when  $z\in N_\delta$ and $n\geq \tilde{N}$. This ensures the
validity of theorem \ref{t4}.

\section{An example of limiting set}

In this section we find the uniform limiting set of second order
spectra for a very simple model. We choose this example so that
the spectral pollution is maximal when we apply the linear method.
We also investigate rank one perturbations of this model. These
ideas are close to example I of \cite{lesh} and those in
\cite{bou}.

Let $E\subsetneqq (-\pi,\pi]$ be a finite union of semi-open
intervals $(a_k,b_k]$. Denote by $E^c:=(-\pi,\pi]\setminus E$. Let
$Mf(x):=m(x)f(x)$ be the linear operator of multiplication by the
symbol
\[
 m(x)=\begin{cases}
    1 & \text{if } x\in E, \\
    -1 & \text{if } x\in E^c
  \end{cases}
\]
acting on $\CH=L^2(-\pi,\pi)$. Let $\CL_n:=\spa \{e^{-i n
x},\ldots,e^{i n x}\}$. Then $M_n$ is the $(2n+1)\times (2n+1)$
Toeplitz matrix $(M_n)_{j,k}=\hat{m}(j-k)$, where $\hat{m}$
denotes the Fourier coefficients of $m$. Notice that $\spec
M=\spec_{\mathrm{ess}} M=\{\pm 1\}$.

Denote by $T(m)$ the Toeplitz operator acting on $l^2(\N)$ whose
symbol is $m$. By invoking Gohberg's theorem for piecewise
continuous symbols (cf., for instance \cite[theorem 1.23]{bosi}),
we realize that $\spec T(m)=[-1,1]$. The truncation of $T(m)$ to
the subspace
\[\tilde{\CL_n}=\spa \{\delta_1,\ldots,\delta_{2n+1}\}\subset
l^2(\N),\] where $\delta_k(n)=\delta_{k,n}$ is the Kronecker delta
symbol, equals $M_n$ above. By virtue of a result due to Szeg\"o
(cf. \cite[theorem 5.14]{bosi}), for each $\mu\in [-1,1]$ there is
a sequence $\mu_n \in \spec M_n$ such that $\mu_n\to \lambda$ as
$n\to \infty$. Thus, if we choose the linear method with $\CL_n$
as approximating subspaces, \emph{each point in the spectral gap
$(-1,1)$ of $M$ is of spectral pollution}.

We now find the limiting set $\Lambda(M)$. According to arguments
in the proof of proposition \ref{t1}, $\spectwo(M|\CL_n)$ lies in
$\T$ for all $n$. We show that each
$\zeta=e^{i\theta},\,-\pi<\theta\leq \pi,$ is in $\Lambda$.
Indeed, by virtue of Gohberg's theorem for piecewise continuous
symbols, the spectrum of $T(m-\zeta)^2$, the Toeplitz operator
associated to the symbol $(m-\zeta)^2$, equals the segment
\[
   \{(1+\zeta)^2t+(1-\zeta)^2(1-t):0\leq t \leq 1\}.
\]
Let $Q_n$ denote the orthogonal projection onto $\tilde{\CL_n}$.
Then
\begin{equation*}
P_n(M-z)^2|\CL_n=Q_nT(m-\zeta)^2|\tilde{\CL}_n.
\end{equation*}
Since the essential spectrum of $T(m-\zeta)^2$ is a set with no
interior, then each point in it is a Weyl's approximated
eigenvalue (cf. \cite[theorem 10.10]{Hisi}). In particular so is
the origin. Hence for all $k>0$, there is $\psi_k \in l^2(\N)$,
$\|\psi_k\|=1$, such that $\|T(m-\zeta)^2\psi_k\|<1/k$. Since
\begin{equation*}
   \sigma_n(\zeta) \leq \frac{\|Q_nT(m-\zeta)^2Q_n\psi_k\|}
   {\|Q_n \psi_k\|}  \to \|T(m-\zeta)^2\psi_k\|
\end{equation*}
as $n\to \infty$. Then $\zeta \in \Sigma$. Consequently
proposition \ref{t1} ensures that $\Lambda(M)=\T.$

\medskip

Let us discuss now how the points of $\spectwo(M|\CL_n)$
distribute along $\T$. We saw in section~2 that
\[
   z\in \spectwo(M|\CL_n) \qquad \text{iff}\qquad \frac{z+z^{-1}}{2}\in
   \spec M_n.
\]
By construction, $M_n$ coincides with the truncation of $T(m)$ to
the subspace $\tilde{\CL}_n$. Let
\[
   \{z_k^{(n)}\}_{k=1}^n:=\spectwo(M|\CL_n)\cap \{\Im z> 0\}
\]
and
\[
   w_k^{(n)}:=\frac{z_k^{(n)}+\overline{z}_k^{(n)}}{2}\in (-1,1).
\]
Then the mean of the points in $\spectwo(M|\CL_n)$ is
\begin{align*}
   \frac{1}{2n}\sum_{k=1}^n z_k^{(n)}+\overline{z}_k^{(n)}&=
   \frac{1}{n}\sum_{k=1}^n w_k^{(n)}
   = \frac{1}{n}\mathrm{Tr}\,M_n \\
   &= \hat{m}(0)= \frac{1}{2\pi}(|E|-|E^c|)
\end{align*}
for all $n$. Here $|\cdot|$ denotes the Lebesgue measure.

A version of Szeg\"o's first limit theorem for Toeplitz operators
\cite[theorem 5.10]{bosi}, permit us to say even more about the
localization of $\spectwo(M|\CL_n)$ in the limit $n\to\infty$.
According to \cite[corollary 5.12]{bosi}, for each Borel set
$B\subset [-1,1]$,
\[
\frac{1}{n}\sum_{w_k^{(n)}\in B}1\to \frac{\big|\{x:m(x)\in
B\}\big|}{2\pi}.
\]
Then for all $0<\varepsilon<1$,
\begin{equation} \label{e9}
\frac{1}{n}\sum_{w_k^{(n)}\in[-1,-1+\varepsilon]}1 \to
\frac{|E^c|}{2\pi} \quad \text{and} \quad
\frac{1}{n}\sum_{w_k^{(n)}\in[1-\varepsilon,1]}1 \to
\frac{|E|}{2\pi}
\end{equation}
as $n\to\infty$. Consequently, although each point in $\T$ is in
the limiting set $\Lambda$, for $n$ large almost
$100\cdot|E|/(2\pi)$ percent of the points in $\spectwo(M|\CL_n)$
cluster near 1 and the other $100\cdot|E^c|/(2\pi)$ percent
cluster near $-1$.

\medskip

\begin{table}
\begin{tabular}{|c|c|c|c|c|}
  \hline
  $\lambda$& $n$ & $\Re z_n -|\Im z_n|$& $\Re z_n +|\Im z_n|$ & $\Re z_n -\lambda$\\
  \hline
   -0.61803398 & 85&  -0.64711164 & -0.59156988 &  0.00130677 \\
 & 120&  -0.64232258 &-0.59559824 &  0.00092642\\
 & 155& -0.63930167 &-0.59820046  &  0.00071708\\
 & 190& -0.63717720 &-0.60006016 &   0.00058469\\
 & 225& -0.63557976 &-0.60147516 &   0.00049347\\
  \hline
  1.61803398 &85&  1.58929960& 1.64481953  &   0.00097441\\
    &120&   1.59398716 & 1.64069975 &   0.00069052\\
   &155 &1.59695260&   1.63804631 &  0.00053453\\
  &190& 1.59904216 &  1.63615390  & 0.00043595\\
   &225& 1.60061565 &  1.63471625  & 0.00036803\\
   \hline
\end{tabular}
\caption{Estimation of the isolated eigenvalues of $M+K$ for
$a=1$, $\psi(x)=1$ and $E=(0,\pi]$. The first column corresponds
to the theoretical value.}
\end{table}

In order to envisage the numerical scope of theorem~\ref{t4}, we
can test it against numerical data. For this we consider rank one
perturbations of the above model. Let $\psi$ be a fixed vector
such that $\|\psi\|=1$, let $a>0$ and let $Kf:=a \langle
f,\psi\rangle \psi$ for all $f\in \CH$. We first compute
explicitly the discrete spectrum of $M+K$.

\begin{lemma} \label{t6}
Let $H$ be self-adjoint and such that $\spec
H=\spec_{\mathrm{ess}}H$.  Then the isolated eigenvalue of finite
multiplicity of $H+K$ are non-degenerate and they are the non-zero
solutions $\lambda$ of
\begin{equation} \label{e10}
   \langle(\lambda-H)^{-1}\psi,\psi\rangle=a^{-1}.
\end{equation}
\end{lemma}
\proof If $(H+K-\lambda)\phi=0$ for $\phi\not=0$ and $\lambda\not
\in \spec H$, then
\[(H-\lambda)\phi+a\langle\phi,\psi\rangle\psi=0.\]
If we normalize by $\langle\phi,\psi\rangle =a^{-1},$ then
$\phi=(\lambda-H)^{-1}\psi$. A substitution in the normalizing
identity yields \eqref{e10}. Conversely suppose that \eqref{e10}
holds. By putting $\phi:=a(\lambda -H)^{-1}\phi$, we achieve
$(H+K-\lambda)\phi =0$.

\medskip

Hence, $\lambda\not=\pm 1$ is an eigenvalue of $M+K$ iff
\begin{equation} \label{e11}
   \frac{\mu_\psi(E)}{\lambda-1}+\frac{\mu_\psi(E^c)}{\lambda+1}=a^{-1},
\end{equation}
where $\mu_\psi(B)=(2\pi)^{-1}\int_B |\psi|^2\,\mathrm{d}x$.

In tables~1 and 2 we show some numerical outputs produced by using
the method of second order spectra applied to $M+K$. We found
$z_n$ by adopting the algorithms for computing second order
spectra available at the internet address mentioned in
\cite{lesh}. This value corresponds to the point in
$\spectwo(M+K|\CL_n)$, closest to the isolated non-degenerate
eigenvalue $\lambda$. We found the theoretical value of $\lambda$
by using \eqref{e11}. In both tables it is remarkable that
although bound \eqref{e2} only provides estimation for the first
digit of $\lambda$ at $n=225$, the actual value of $\Re z_n$ is
correct up to three digits. For table~2 we chose an extreme case
where one of the non-degenerate eigenvalues is close to the point
-1 of the essential spectrum. It is remarkable that the even the
step $n=85$, which only takes  a few second to run in a PC,
already detects the presence of this eigenvalue.

\begin{table}
\begin{tabular}{|c|c|c|c|c|}
  \hline
  $\lambda$& $n$ & $\Re z_n -|\Im z_n|$& $\Re z_n +|\Im z_n|$ & $\Re z_n -\lambda$\\
  \hline
  -0.97901994 & 85& -0.99169545 &-0.97384630 &    0.00375093 \\
  & 120& -0.98897219  & -0.97406728&   0.00249979 \\
  &155&-0.98740174 &-0.97435952  &  0.00186068 \\
  &190 &-0.98635662&-0.97465104 &    0.00148388 \\
  &225& -0.98561863&-0.97491625  &  0.00124750\\
  \hline
  1.97901994 &85&1.95326913  & 2.00377483 &  0.00049795\\
   &120&1.95700470 &  2.00030477  & 0.00036520 \\
   &155&1.95961155 &  1.99784090 &  0.00029371\\
   &190&1.96164380  & 1.99591836 &  0.00023886\\
  &225& 1.96314714 &  1.99449943 &  0.00019665\\
   \hline
\end{tabular}
\caption{Estimation of the isolated eigenvalues of $M+K$ for
$a=1$, $\psi(x)=1$ and $E=(-15 \pi/16,\pi]$. The first column
corresponds to the theoretical value.}
\end{table}

Figures~1, 2 and 3, shows how the rest of the second order
spectrum in these two cases distribute along the complex plane for
three different values of $n$. Notice that the clustering
predicted for the unperturbed case seems to be largely kept with
the exception of some few points that approximate the isolated
eigenvalues.

\section{The condition (H)}

As in the linear case, an effective numerical implementation of
the method of second order spectra to unbounded operators, demands
verifying (H) without having much information about the
eigenfunctions $\psi$. This is achieved sometimes by means of a
perturbative argumentation. We do not claim that this is easy in
general, but in some situations the following standard argument
can be useful.

The general strategy involves a dominating operator $X$, which has
compact resolvent and a complete set of eigenfunctions that can be
found explicitly. Both operators $M^k$, $k=1,2$, should be
relatively $X$-bounded, in the sense that $\dom X \subseteq \dom
M^2$ and
\[
   \|M^k u \| \leq a_k\|X u\|+b_k\|u\|, \qquad \qquad
   u\in \dom X,
\]
for uniform constants $a_k>0$ and $b_k>0$. The best $a_k$ such
that the above holds for some $b_k$, is called relative bound. Put
as approximating sequence $\{\phi_n\}_{n=1}^\infty$, the set of
eigenfunctions of $X$. If the eigenfunction  $\psi$ of $M$
associated to the eigenvalue $\lambda$ lies in $\dom X$, by virtue
of the spectral theorem, $X$ and $P_n$ commute, and hence
\[
  \|XP_n\psi -X \psi\|= \|P_n X\psi -X\psi \| \to 0
\]
as $n\to \infty$. Then $\|M^kP_n \psi- M^k\psi\|\to 0$, ensuring
(H).

The following result illustrates this strategy. Recall that the
class $K_1$ consists of all $V:\R\longrightarrow \R$ such that
\[
    \sup_{x\in \R} \int _{|x-y|\leq 1} |V(y)| \mathrm{d} y
<\infty.
\]

\begin{theorem} \label{t8}
Let $\phi_n$ be the eigenfunctions of the harmonic oscillator
$(-\partial_x^2+x^2)$ acting on $L^2(\R)$. Let
$V:\R\longrightarrow \R$ be such that \[V\in [K_1\cap
W^{2,2}(\R)]+W^{\infty,2}(\R).\] Then (H) holds for every
eigenvalue and eigenfunction of the operator $M=-\partial_x^2+V$
acting on $L^2(\R)$.
\end{theorem}

Analogous results can be found for higher dimensions.

\proof  Since $V\in L^2+L^\infty$, multiplication by $V$ is
relatively $\partial^2_x$-bounded with bound 0 (cf. \cite[\S
1.2]{cfks}), hence $M$ is relatively $\partial_x^4$-bounded with
bound 0.

We prove that $M^2$ is relatively $\partial_x^4$-bounded with
bound 1. For this, notice that
\[
   (-\partial_x^2+V)^2u =
   \partial_x^4u-2V\partial_x^2u-2V'\partial_xu+(V^2-V'')u.
\]
Since $V\in L^2+L^\infty$, multiplication by $V$ is
$\partial_x^2$-bounded with relative bound 0 and hence
$2V\partial_x^2$ is $\partial_x^4$-bounded with relative bound 0.
Since $V'\in L^2+L^\infty$, a similar reasoning ensures that
$2V'\partial_x$ is $\partial_x^3$-bounded with relative bound 0.
Since $V''$ lies in $L^2+L^\infty$, multiplication by $V''$ is
$\partial_x^2$-bounded with relative bound 0. Since $V^2$ is
$\partial_x^2V$-bounded with relative bound 0 and
\[
  \partial_x^2Vu=V''u+2V'\partial_xu+V\partial_x^2u,
\]
then $V^2$ is $\partial_x^4$-bounded with relative bound 0. This
ensures that $M^2$ is $\partial_x^4$-bounded with relative bound
1.

By virtue of Kato-Rellich theorem, above we do not have to worry
about the domains. Here we are considering
\[
   \dom M=W^{2,2} \qquad \mathrm{and} \qquad
   \dom M^2=W^{2,4}.
\]
The mentioned result ensures that $M$ (and hence $M^2$) are
self-adjoint in these domains.

Put $X=(-\partial_x^2+x^2)^2$. Since $V\in K_1$, then
$|\psi(x)|\leq c e^{-a|x|}$ for some constants $a>0$ and $c>0$
(cf. \cite[\S C.3]{250}). Together with the inclusion $\psi\in
W^{2,4}$, this bound ensures $\psi\in \dom
(-\partial_x^2+x^2)=W^{2,2}\cap \dom(x^2)$ and
$(-\psi''+x^2\psi)\in W^{2,2}\cap \dom(x^2)$. The only non-trivial
facts of the latter assertion are, perhaps, the inclusions
$x^2\psi\in W^{2,2}$ and $\psi''\in \dom(x^2)$. The first follows
from the second, by differentiating twice the term $x^2\psi$ and
noticing that $(x\psi')'\in L^2$. The second inclusion is achieved
by means of the following trick. Since $M\psi$ is eigenvector of
$M$,
$
   |-\psi''(x)+V(x)\psi(x)|\leq c e^{-a|x|}.
$
Then
\begin{align*}
  \|x^2 \psi''\|& \leq \|x^2(-\psi''+V\psi)\|+\|x^2V\psi\| \\
   & \leq c_1+ \left(\int x^4 |V(x)|^2\,|\psi(x)|^2
   \mathrm{d} x\right)^{1/2} \\
   & \leq c_1 + c_2\left( \int x^4 e^{-a|x|} |V(x)|^2 \mathrm{d} x\right)^{1/2}.
\end{align*}
The latter integral is bounded because of $V\in L^2+L^\infty$,
therefore \linebreak $\psi''\in \dom(x^2)$. Thus $\psi\in \dom X$.

Finally we show that $\partial_x^4$ is
$(\partial_x^2+x^2)^2$-bounded. For this, let
\[
    A:=2^{-1/2}(x+\partial_x)\qquad \mathrm{and} \qquad
    A^\ast = 2^{-1/2}(x-\partial_x).
\]
Then $(-\partial_x^2+x^2)=2(AA^\ast-1)$ and
$\partial_x^4=(A-A^\ast)^4$. Thus the desired property follows
from the identity (cf.\cite[eq.(X.28)]{resiv2})
\[
    \|A^{\#_1}\cdots A^{\#_n}u\|\leq c
    \|(-\partial_x^2+x^2)^{n/2}u\|, \qquad n=1,2,\ldots
\]
where $A^{\#_k}$ is either $A$ or $A^\ast$. This is
easily shown by induction and using
the estimate
\[
  \|(-\partial_x^2+x^2)^{k/2}u\|\leq\|(-\partial_x^2+x^2)^{n/2}u\|,
  \qquad k<n.
\]
This completes the proof of the theorem.

\medskip

In the above result, we have in mind the case $V=W+S$, where $W$
is periodic and $S$ is in $L^2(\R)+L^\infty(\R)_\varepsilon$.
Under this hypothesis, it is well known  that the essential
spectrum of $M$ has a band gap structure determined solely by $W$,
whereas the perturbation $S$ can produce  non-empty discrete
spectrum (cf. e.g. \cite[Chapter XIII.16]{resi} and also \cite{ki}).


\bigskip

{\samepage {\scshape Acknowledgments.} The author wishes to thank
E.B.~Davies, E.~Shargorodsky and V.~Strauss for their useful comments.}

\vspace{1in}

{\scshape \centerline{Departamento de Matem\'aticas Puras y
Aplicadas,}} {\scshape \centerline{Universidad Sim\'on Bol\'\i
var,}}  {\scshape \centerline{Apartado 89000, Caracas 1080-A,
Venezuela.}}

\centerline {\textit{E-mail address:} \texttt{lboulton@ma.usb.ve}}

\clearpage

\begin{figure}[t] \label{f1}
\begin{picture}(250,250)(75,130) \includegraphics{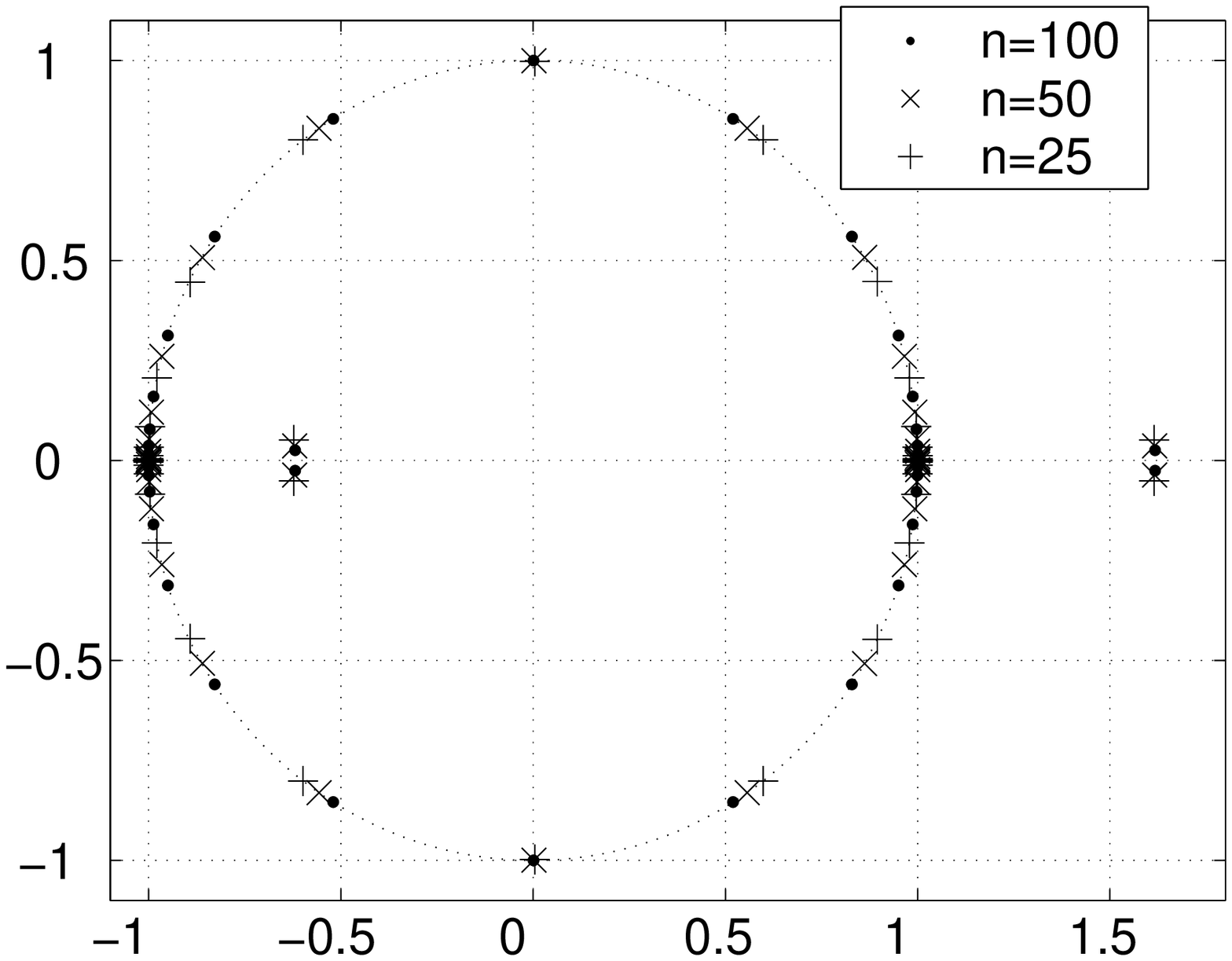}
\end{picture}
\caption{$\spectwo(M+K|\CL_n)$ for three values of $n$. The
operator corresponds to $a=1$, $\psi(x)=1$ and $E=(0,\pi]$. See
table~1. }
\end{figure}

\begin{figure}[t] \label{f2}
\begin{picture}(250,250)(75,130) \includegraphics{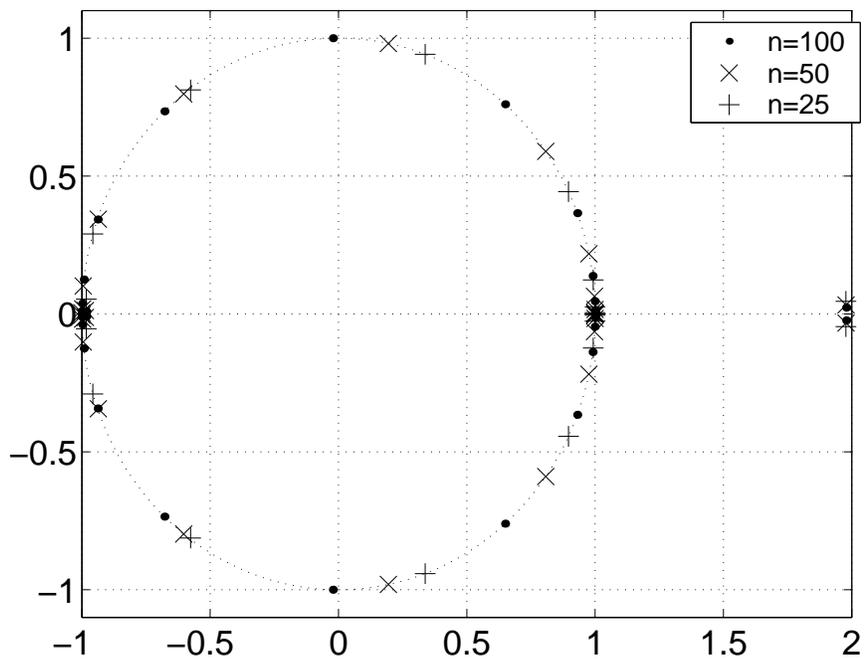}
\end{picture}
\caption{$\spectwo(M+K|\CL_n)$ for three values of $n$. The
operator corresponds to $a=1$, $\psi(x)=1$ and $E=(-15\pi/16,
\pi]$. See table~2.}
\end{figure}

\begin{figure}[b] \label{f3}
\begin{picture}(250,200)(75,200) \includegraphics{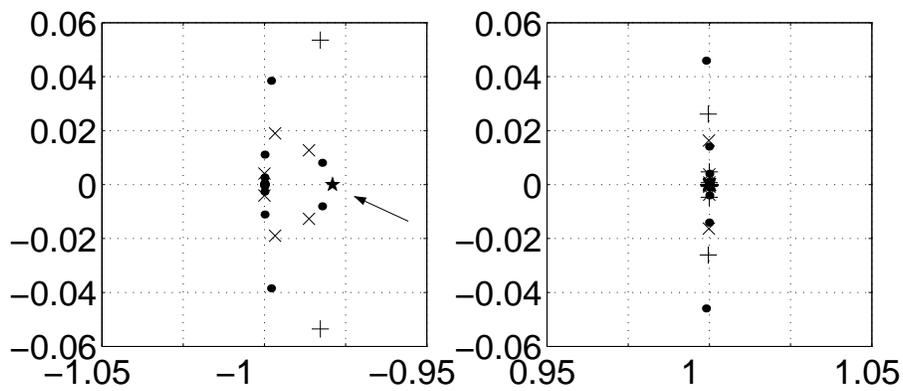}
\end{picture}
\caption{Two zoom pictures of figure~2. Clustering near $\pm 1$.
The star shows the location of the theoretical eigenvalue.}
\end{figure}

\end{document}